\documentclass[12pt,a4paper]{article}
\usepackage{amscd}
\usepackage[latin1]{inputenc}
\usepackage{amsmath}
\usepackage{amsfonts}
\usepackage{amssymb}
\usepackage{color}
\usepackage{float}
\usepackage{amsthm}

\newtheorem{theorem}{Theorem}

\newtheorem{lemma}[theorem]{Lemma}
\newtheorem{proposition}[theorem]{Proposition}
\newtheorem{remark}[theorem]{Remark}

\newtheorem*{minmod}{(*)}

\def\Qed{\hfill\raisebox{.6ex}{\framebox[2.5mm]{}}\\[.15in]}
\def\b{\overline}
\def\m{\mathbb}
\def\mc{\mathcal}

\begin{document}
\date{}
\title{Hyperelliptic surfaces with $K^2 < 4\chi - 6$}
\author{Carlos Rito \and Mar\'ia Mart\'i S\'anchez}
\maketitle

\begin{abstract}

Let $S$ be a smooth minimal surface of general type with a (rational) pencil of hyperelliptic curves of minimal genus $g$.
We prove that if $K_S^2<4\chi(\mathcal O_S)-6,$ then $g$ is bounded.
The surface $S$ is determined by the branch locus of the covering $S\rightarrow S/i,$ where $i$ is the hyperelliptic involution of $S.$
For $K_S^2<3\chi(\mathcal O_S)-6,$ we show how to determine the possibilities for this branch curve.
As an application, given $g>4$ and $K_S^2-3\chi(\mathcal O_S)<-6,$ we compute the maximum value for $\chi(\mathcal O_S).$ This list of possibilities is sharp.

\end{abstract}

\section{Introduction}

Our motivation for this work is the following. In \cite{AK} Ashikaga and Konno consider surfaces $S$ of general type with $K_S^2=3\chi(\mathcal O_S)-10.$
For these surfaces the canonical map is of degree 1 or 2. In the degree 2 case, the canonical image is a ruled surface, thus if $S$ is regular, it has a pencil of hyperelliptic curves.
By a result of Xiao \cite[Thm. 1]{Xi2} if $\chi(\mathcal O_S)\geq 47,$ then $S$ has such an hyperelliptic pencil of curves of genus $\leq 4.$
But for $\chi(\mathcal O_S)\leq 46$ this result gives no information
(for $\chi(\mathcal O_S)=46$ the slope formula \cite[Thm. 2]{Xi1} implies $g\leq 5\vee g\geq 9;$  we show that in this case $S$ has an hyperelliptic pencil of minimal genus $g\leq 10$ and the cases $g=9,$ $g=10$ do occur).
Ashikaga and Konno study only the case $g\leq 4$ (there is an infinite number of possibilities). Nothing is said for the possibilities with $g\geq 5$ and $\chi(\mathcal O_S)\leq 46.$
A similar situation occurs in \cite{K}.

In this paper we study smooth minimal surfaces $S$ of general type which have a pencil of hyperelliptic curves (by {\em pencil} we mean a linear system of dimension 1).
We say that $S$ has such a pencil of {\em minimal genus} $g$ if it has an hyperelliptic pencil of genus $g$ and all hyperelliptic pencils of $S$ are of genus $\geq g.$ We are mainly interested in the case $g>4$ and $\chi(\mathcal O_S)$ small (i.e. where \cite[Thm. 1]{Xi2} is not useful).

For $S$ such that $K_S^2<4\chi(\mathcal O_S)-6,$ we give bounds for the minimal genus $g$ (Theorem \ref{thm1}).

The surface $S$ is the smooth minimal model of a double cover of an Hirzebruch surface $\m F_e$ ramified over a curve $\b B$ (which determines $S$).
We prove that if $K_S^2<3\chi(\mathcal O_S)-6,$ then $\b B$ has at most points of multiplicity $8$ and we show how to determine the possibilities for $\b B$ (Proposition \ref{thm2}).

As an application, given $g>4$ and $K_S^2-3\chi(\mathcal O_S)<-6,$ we compute the maximum value for $\chi(\mathcal O_S);$ this list of possibilities is sharp (Theorem \ref{thm3}).

The paper is organized as follows.
In Section \ref{MR} we present the main results of the paper.
The hyperelliptic involutions of the fibres of $S$ induce an involution $i$ of $S$, so in Section \ref{GenFacts} we review some general facts on involutions.
Since the quotient $S/i$ is a rational surface, a smooth minimal model of $S/i$ is not unique.
We make a choice for this minimal model in Section \ref{choice} (which is due to Xiao \cite{Xi3}) and we show some consequences of it.
Section \ref{BG} contains the key result of the paper, which allow us to compute bounds for the minimal genus of the hyperelliptic fibration.
This is done via a carefully analysis of the possibilities for the branch locus of the covering $S\rightarrow S/i$ considering the restrictions imposed by the choice of minimal model.
Finally this is used in Section \ref{PT} to prove the main results, stated in Section \ref{MR}.

\bigskip
\noindent{\bf Notation}

We work over the complex numbers; all varieties are assumed to be projective algebraic.
A {\em $(-2)$-curve} or {\em nodal curve} $A$ on a surface is a curve isomorphic to $\mathbb P^1$ such that $A^2=-2$.
An $(m_1,m_2,\ldots)$-point of a curve, or point of type $(m_1,m_2,\ldots),$ is a singular point of multiplicity $m_1,$ which resolves to a point of multiplicity $m_2$ after one blow-up, etc.
By {\em double cover} we mean a finite morphism of degree $2.$
The rest of the notation is standard in Algebraic Geometry.\\

\bigskip
\noindent{\bf Acknowledgements}

The authors wish to thank Margarida Mendes Lopes, for all the support, Ana Bravo and Orlando Villamayor, for their hospitality at the Universidad Aut\'onoma de Madrid.

The first author is a member of the Mathematics Center of the Universidade de Tr\'as-os-Montes e Alto Douro. Both authors are collaborators of the Center for Mathematical Analysis, Geometry and Dynamical Systems of Instituto Superior T\'ecnico, Universidade T\' ecnica de Lisboa.

This research was partially supported by the Funda\c c\~ao para a Ci\^encia e a Tecnologia (Portugal) through Projects PEst-OE/MAT/UI4080/2011 and PTDC/MAT/099275/2008.

\section{Main results}\label{MR}

\begin{theorem}\label{thm1}

Let $S$ be a minimal smooth surface of general type with a pencil of hyperelliptic curves of minimal genus $g$.

If $K_S^2 < 4\chi(\mathcal O_S)-6,$
then $g$ is not greater than
\\
$\max \left\{ -1+\frac{8\chi(\mathcal O_S)}{4\chi(\mathcal O_S)-K_S^2-6},1+\frac{8\chi(\mathcal O_S)-16}{4\chi(\mathcal O_S)-K_S^2-6} ,1+\frac{8\chi(\mathcal O_S)}{4\chi(\mathcal O_S)-K_S^2-3}, \frac{3+\sqrt{1+8\chi(\mathcal O_S)}}{2} \right\}.$

\end{theorem}

Let $B\subset W$ be the branch locus of a double cover $V\rightarrow W,$ where $V$ and $W$ are smooth surfaces (thus $B$ is also smooth). Let $\rho:W\rightarrow P$ be the projection of $W$ onto a minimal model and denote by $\overline B$ the projection $\rho(B).$

Suppose that $\overline B$ has singular points $x_1,\ldots,x_n$ (possibly infinitely near). For each $x_i$ there is an exceptional divisor $E_i$ and a number $r_i\in 2\mathbb N$ such that
$$
\begin{array}{l}
E_i^2=-1,\\
K_W\equiv\rho^*(K_P)+\sum E_i,\\
B=\rho^*(\overline{B})-\sum r_iE_i .
\end{array}
$$

Notice that $r_i$ is not the multiplicity of the singular point $x_i,$ it is the multiplicity of the corresponding singularity in the {\em canonical resolution} (see \cite[III. 7.]{BHPV}).
For example, in the case of a point of type $(2r-1,2r-1)$ one has $r_1=2r-2$ and $r_2=2r.$

Since, from Theorem \ref{thm1}, we have a bound for the genus $g,$ we also have a bound for the multiplicities $r_i.$ For the case $K_S^2<3\chi(\mathcal O_S)-6,$ we prove the result below.

Let $N_j$ be the number of singular points $x_i$ of $\b B$ (possibly infinitely near) such that $r_i=j.$
Denote by $C_0$ and $F$ the negative section and a ruling of the Hirzebruch surface $\mathbb F_e.$

\begin{proposition}\label{thm2}
Let $S$ be a minimal smooth surface of general type with an hyperelliptic pencil of minimal genus $(k-2)/2$.
If $K_S^2<3\chi(\mathcal O_S)-6,$ then $S$ is the smooth minimal model of a double cover $S'\rightarrow\mathbb F_e$ with branch curve $\b B\equiv kC_0+(ek/2+l)F$ such that:
\begin{description}
\item[a)] $r_i\leq\min\{8,\ k/2+2,\ l-k/2+2\}$\ $\forall i$
\item[b)] $N_4+N_6=15+K_{S''}^2-3\chi(\mathcal O_S)-\frac{1}{4}(k-10)(l-10)$
\item[c)] $\chi(\mathcal O_S)=1+\frac{1}{4}(k-2)(l-2)-N_4-3N_6-6N_8$
\end{description}
where $S''\longrightarrow S'$ is the canonical resolution.
\end{proposition}

Proposition \ref{thm2} can be used to restrict possibilities for $\b B.$ We show the following:

\begin{theorem}\label{thm3}
Let $S$ be a minimal smooth surface of general type with an hyperelliptic pencil of minimal genus $g>4$.
If $K_S^2 < 3\chi(\mathcal O_S) - 6,$ then $\chi(\mathcal O_S)$ is bounded by the number given in the table below (emptiness means non-existence).
All these cases do exist.

\begin{table}[H]\label{tablechi}
\centering
\begin{tabular}{l l||l|l|l|l|l|l|l|l|l|l}
        ~                & g & ~   & ~   & ~   & ~ & ~ & ~ & ~  & ~ & ~ \\
        $K^2-3\chi$ & ~ &      -7  & -8  & -9  & -10  & -11  & -12  & -13  & -14  & -15 & -16 \\\hline\hline
        ~                & 5 & 61  & 56  & 51  &  46  & 41   & 36   & 31   & 26   & 21  & 16  \\\hline
        ~                & 6 & 49  & 46  & 43  &  40  & 37   & 34   & 27   & 28   & ~   & 22  \\\hline
        ~                & 7 & 42  & 43  & 43  &  35  & 35   & 36   & 28   & ~    & 29  & 22  \\\hline
        ~                & 8 & 44  & 44  & 45  &  ~   & 36   & ~    & 37   & ~    & 29  & ~   \\\hline
        ~                & 9 & ~   & 45  & ~   &  46  & ~    & ~    & 37   & ~    & ~   & ~   \\\hline
        ~               & 10 & ~   & ~   & ~   &  46  & ~    & ~    & ~    & ~    & ~   & ~   \\     
\end{tabular}
\end{table}
\end{theorem}

\begin{remark}
This result gives $3$ examples where Theorem \ref{thm1} is almost sharp:
in the cases $\left(g,K^2-3\chi\right)=(10,-10), (9,-13), (8,-15)$ we get from Theorem \ref{thm1} that $\chi\leq 47, 38, 30,$ respectively.

There is at least one case where Theorem \ref{thm1} is sharp: a double plane with branch locus a curve of degree $18$ with $8$ points of multiplicity $6.$
In this case $\chi=5,$ $K^2=8$ and $g=5.$
\end{remark}

\section{Involutions}\label{GenFacts}

Let $S$ be a smooth minimal surface of general type with a (rational) pencil of hyperelliptic curves.
This hyperelliptic structure induces an involution (i.e. an automorphism of order 2) $i$ of $S$.
The quotient $S/i$ is a rational surface.

Since S is minimal of general type, this involution is biregular. 
The fixed locus of $i$ is the union of a smooth curve $R''$ (possibly empty) and of $t\geq 0$ isolated points $P_1,\ldots,P_t.$ Let $p:S\rightarrow S/i$ be the projection onto the quotient. The surface $S/i$ has nodes at the points $Q_i:=p(P_i),$ $i=1,\ldots,t,$ and is smooth elsewhere. If $R''\not=\emptyset,$ the image via $p$ of $R''$ is a smooth curve $B''$ not containing the singular points $Q_i,$ $i=1,\ldots,t.$ Let now $h:V\rightarrow S$ be the blow-up of $S$ at $P_1,\ldots,P_t$ and set $R'=h^*(R'').$ The involution $i$ induces a biregular involution $\widetilde{i}$ on $V$ whose fixed locus is $R:=R'+\sum_1^t h^{-1} (P_i).$ The quotient $W:=V/\widetilde{i}$ is smooth and one has a commutative diagram
$$
\begin{CD}\ V@>h>>S\\ @V\pi VV  @VV p V\\ W@>g >> S/i
\end{CD}
$$
where $\pi:V\rightarrow W$ is the projection onto the quotient and $g:W\rightarrow S/i$ is the minimal desingularization map. Notice that $$A_i:=g^{-1}(Q_i),\ \ i=1,\ldots,t,$$ are $(-2)$-curves and $\pi^*(A_i)=2\cdot h^{-1}(P_i).$ 

Set $B':=g^*(B'').$ Since $\pi$ is a double cover, its branch locus $B'+\sum_1^t A_i$ is {\em even}, i.e. there is a line bundle $L$ on $W$ such that $$2L\equiv B:=B'+\sum_1^t A_i.$$

\section{Choice of minimal model}\label{choice}

Part of this section may be found in \cite{Xi3}. We use the notation introduced so far.
As above, $W$ is a rational surface.
\begin{minmod}
Blowing-up, if necessary, $\mathbb P^2$ at a point, we can suppose that $W\ne\mathbb P^2.$
\end{minmod}
\noindent Thus there is a birational morphism $$\rho:W\longrightarrow\mathbb F_e,$$ where $\mathbb F_e$ is an Hirzebruch surface.
Let $\overline B:=\rho(B)$ and consider the double cover $S'\longrightarrow\mathbb F_e$ with branch locus $\overline B.$ If $\overline B$ is singular then $S'$ is also singular and $S$ is isomorphic to the minimal smooth resolution of $S'$. 

We can define $k$ and $l$ such that $$\overline{B}\equiv :kC_0+\left(\frac{ek}{2}+l\right) F,$$
where $C_0$ and $F$ are, respectively, the negative section and a ruling of $\mathbb F_e$ (thus $C_0^2=-e,$ $C_0F=1,$ $F^2=0$).
Notice that $\b B^2=2kl$ and $K_P\b B=-2k-2l$.

\begin{minmod}
Among all the possibilities for the map $\rho,$ we choose one satisfying, in this order:
\begin{description}
  \item[1)]  the degree $k$ of $\overline{B}$ over a section is
minimal;
  \item[2)]  the greatest order of the singularities of $\overline{B}$
is minimal;
  \item[3)]  the number of singularities with greatest order is
also minimal.
\end{description}
\end{minmod}

Recall that a $(2r-1, 2r-1)$ singularity of $\overline{B}$ is a pair $(x_j,x_k)$ such that $x_k$ is infinitely near to $x_j$ and $r_j=2r-2,$ $r_k=2r$. 

Let $$r_m:=\max\ \{r_i\}$$ or $r_m:=0$ if $\b B$ is smooth.

By {\em elementary transformation} over $x_i\in\m F_e$ we mean the blow-up of $x_i$ followed by the blow-down of the strict transform of the ruling of $\m F_e$ that contains $x_i.$

The following is a consequence of the choice (*) of the map $\rho$.

\begin{proposition}[\cite{Xi3}]\label{semiXG}
We have:
\begin{description}
\item[a)] If $k\equiv 0\ ({\rm mod\ 4}),$ then $r_m\leq \frac{k}{2}+2$ and the equality holds only if $x_m$ belongs to a singularity
$\left(\frac{k}{2}+1,\frac{k}{2}+1\right).$ In this case $l\geq k+2$ and all the branches of the singularity are tangent to the ruling of $\m F_e$ that contains it.
\item[b)]  If $k\equiv 2\ ({\rm mod\ 4}),$ then $r_m\leq \frac{k}{2}+1$ and the equality holds only if $x_m$ belongs to a singularity
$\left(\frac{k}{2},\frac{k}{2}\right).$ In this case $l\geq k.$
\end{description}
\end{proposition}

In a similar vein:

\begin{proposition}\label{semiXG2}
We have that:

\begin{description}
\item[a)] if $l=k+2$ and $k>8,$ there are at most two $\left(\frac{k}{2}+1,\frac{k}{2}+1\right)$-points.
\item[b)] $l\geq\frac{k}{2}$ and $l\geq\frac{k}{2}+r_m-2;$
\item[c)] if $l=\frac{k}{2}+r_m-2,$ then either:
\begin{description}
\item[$\cdot$] $e=2,$ $l=k-2,$ the branch locus $\b B$ has a $\left(\frac{k}{2}-1,\frac{k}{2}-1\right)$-point and all singularities are of multiplicity $<\frac{k}{2},$ or
\item[$\cdot$] we can suppose $e=1,$ the negative section $C_0$ of $\m F_1$ is contained in $\b B,$ $\b B$ has a point of multiplicity $r_m$ contained in $C_0$ and the remaining singularities are of multiplicity $<r_m.$
\end{description}
\end{description}

\end{proposition}

\noindent{\bf Proof:}
\begin{description}
\item[a)] This is due to Borrelli (\cite{Bo}). 
Suppose that there are three singularities $\left(k/2+1,k/2+1\right).$
The rulings of $\m F_e$ through these points are contained in $\b B$ and then $\b BC_0=l-\frac{ek}{2}\geq 4$ ($\b BC_0$ is even). This implies $e\leq 1.$ Making, if necessary, an elementary transformation over one of these points, we can suppose that $e=1.$

Let $\rho$ be as above and $E_i, E_i',$ $i=1,2,3,$ be the exceptional divisors corresponding to three singularities $(k/2+1,k/2+1)$ of $\b B.$
The general element of the linear system ${|\rho^*(4C_0+5F)-\sum_1^3(2E_i+2E_i')|}$ is a smooth and irreducible rational curve $C$ such that $CB<k.$
This contradicts the choice (*) of the map $\rho.$

\item[b)] If $r_m>\frac{k}{2}$ then the result follows from Proposition \ref{semiXG}. Suppose
now $r_m\leq \frac{k}{2}.$ We have $\overline{B}C_0\geq -e,$ i.e.
$l-\frac{ek}{2}\geq -e.$
Therefore if $e\geq 2,$ then $$l\geq k-2\geq\frac{k}{2} {\rm\ \ \ \ and\ \ \ \ } l\geq k-2\geq\frac{k}{2}+r_m-2.$$

When $e=0$ we obtain immediately $l\geq k,$ by the choice of the map $\rho,$ thus $l\geq\frac{k}{2}+r_m.$

If $e=1$ then $\overline{B}C_0=l-\frac{k}{2}\geq 0.$ Blow-down $C_0.$ We
obtain a singularity of order at most $l-\frac{k}{2}+1,$ hence the
choice of the minimal model implies $r_m\leq l-\frac{k}{2}+2$ (notice
that the equality happens only if the order of the singularity is
$\left(r_m-1,r_m-1\right)$). 

\item[c)] 
Assume that $l=k/2+r_m-2.$ Proposition \ref{semiXG} implies $r_m\leq k/2.$
From $\b BC_0\geq -e$ we obtain $k/2+r_m-2=l\geq\frac{ek}{2}-e,$ thus either $e=1$ or $e=2$ and $r_m=k/2$
(notice that $e=0$ implies $l\geq k$).

In the case $e=1$ we can, as in the proof of b), contract the section with self-intersection $(-1)$ to obtain a branch curve in $\mathbb P^2$ with at most singularities of type
$\left(l-k/2+1,l-k/2+1\right).$

Suppose now that $e=2$ and there is a point $x_i$ of multiplicity $k/2.$
In this case $\b BC_0=-2,$ hence $x_i\not\in C_0.$ We make an elementary transformation over $x_i$ to obtain the case $e=1$ also with $l=k-2.$ \Qed

\end{description}

\section{Bound of genus}\label{BG}

In this section we prove the key result to establish bounds for the minimal genus of the hyperelliptic fibrations.

From \cite{Ri} (cf. also \cite{CM}), we get the following:
\begin{proposition}\label{rito}

Let $S''\rightarrow S'$ be the canonical resolution of a double cover $S'\rightarrow\m F_e$ with branch locus $\b B\equiv kC_0+(ek/2+l)F$. Let $S$ be the minimal model of $S''$ and $t:=K_S^2-K_{S''}^2.$ If $S$ is of general type, then:
\
\begin{description}
\item[a)] $\sum (r_i-2)(k-r_i-2)=H$
\item[b)] $2l=G+\sum(r_i-2),$
\end{description}
where\\
$H=2k^2-k\left( 4\chi(\mathcal O_S)+t-K_S^2+8 \right)+16\chi(\mathcal O_S)+2t-2K_S^2$\\
and
\\
$G=-2k+4\chi(\mathcal O_S)+t-K_S^2+8$.
\end{proposition}
\noindent{\bf Proof:} 
From \cite[Propositions 2 and 3, a)]{Ri} one gets:
\begin{description}
\item[(a)] $2kl=-48+12l+12k-8\chi(\mathcal O_S)+4K_S^2-4t+ \sum (r_i-2)(r_i-4)$
\item[(b)] $2k+2l=8+4\chi(\mathcal O_S)+t-K_S^2 + \sum(r_i-2).$
\end{description}
The result is obtained replacing (a) by (a)+$(6-k)$(b). \Qed

The next result is a fundamental tool in the proof of Proposition \ref{mainprop} below.

\begin{lemma}\label{mainlemma}
Suppose that $k> 8$. With the above notation, we have
\
\begin{description}
\item[a)] $2l\leq G+\dfrac{H}{k-r_m-2},$ and
\item[b)] if $r_m$ is obtained only from singularities of type $(r_m -1, r_m -1)$, then
$$2l\leq G+\dfrac{H}{(r_m-4)(k-r_m)+(r_m-2)(k-r_m-2)}(2r_m-6).$$
\end{description}
\end{lemma}
\noindent{\bf Proof:}
The first statement follows from Proposition \ref{rito} and Proposition \ref{semiXG}, a). 

Next we prove b). By the assumptions,  if $x_i$ does not belong to a $(r_m -1, r_m -1)$ singularity, we have $r_i< r_m$.   
Let $n\geq 1$ be the number of singularities of type $(r_m-1, r_m-1)$ and $s\geq 0$ be the number of singular points $x_j$ of another type.
As seen in Section 4, each singularity $(r_m-1, r_m-1)$ corresponds to two infinitely near singular points $x_k,$ $x_{k+1}$ with $r_k=r_m-2,$ $r_{k+1}=r_m$.
Therefore $$\sum_{i=1}^{2n+s} (r_i-2)=n(2r_m-6) + \sum_{j=1}^{s}(r_j-2),$$ with $r_j<r_m.$
Thus from Proposition \ref{rito}, b) we get
$$2l=G + n(2r_m-6)+ \sum_{j=1}^{s}(r_j-2).$$

By Proposition \ref{rito}, a),
$$H=n\Big((r_m-4)(k-r_m) + (r_m-2)(k-r_m-2)\Big) + \sum_{j=1}^{s}(r_j-2)(k-r_j-2),$$
hence $$n=\frac{H-\sum_{j=1}^{s}(r_j-2)(k-r_j-2)}{(r_m-4)(k-r_m) + (r_m-2)(k-r_m-2)}$$
and then
\begin{equation}\label{eqmaria}
2l=G + \frac{H-\sum_{j=1}^{s}(r_j-2)(k-r_j-2)}{(r_m-4)(k-r_m) + (r_m-2)(k-r_m-2)}(2r_m-6)+\sum_{j=1}^{s}(r_j-2).
\end{equation}

Since $r_j< r_m,$ $j=1,\ldots,s,$ 
$$(r_m-4)(k-r_m) + (r_m-2)(k-r_m-2) \leq (2r_m-6)(k-r_j-2).$$ This implies 
$$\sum_{j=1}^{s}(r_j-2)\leq \sum_{j=1}^{s}\frac{(r_j-2)(k-r_j-2)(2r_m-6)}{(r_m-4)(k-r_m) + (r_m-2)(k-r_m-2) }$$
and the result follows from (\ref{eqmaria}). \Qed

The following proposition will allow us to give bounds for $k.$
Notice that, since $\b B$ is even and $\b BC_0=l-\frac{ek}{2},$ $$k\equiv 0\ ({\rm mod}\ 4)\ \Longrightarrow\ l\equiv 0\ ({\rm mod}\ 2).$$
 
\begin{proposition}\label{mainprop}
In the conditions of Proposition \ref{rito}, suppose that $k>8.$
\begin{description}
\item If $k\equiv 0\ ({\rm mod\ 4}),$ one of the following holds:
\begin{description}
\item[a)] $r_m=k/2+2,$ $l=k+2$ and\newline
$\left( 4\chi(\mathcal O_S)+t-K_S^2-8 \right)k\leq 16\chi(\mathcal O_S) - 16,$ with $t\geq 2$;
\item[b)] $r_m=k/2+2,$ $l\geq k+4$ and\newline
$(4\chi(\mathcal O_S) + t - K_S^2 - 8)k^2-16\chi(\mathcal O_S)k+32\chi(\mathcal O_S)\leq 0,$ with $t\geq 2$;
\item[c)] $r_m=k/2,$ $l=k-2$ and\newline
$(4\chi(\mathcal O_S) + t - K_S^2 - 4)k^2+(-48\chi(\mathcal O_S) - 8t + 8K_S^2 + 32)k+$\newline $160\chi(\mathcal O_S) + 16t - 16K_S^2 - 96\leq 0,$ with $t\geq 1$,\newline
or\newline
$(4\chi(\mathcal O_S) + t - K_S^2 + 2)k\leq 32\chi(\mathcal O_S) + 4t - 4K_S^2 - 8,$ with $t\geq 1$,\newline
or\newline
$(4\chi(\mathcal O_S) + t - K_S^2 - 5)k^2+(-48\chi(\mathcal O_S) - 8t + 8K_S^2 + 44)k+$\newline $160\chi(\mathcal O_S) + 16t - 16K_S^2 - 128\leq 0,$ with $t\geq 2$;
\item[d)] $r_m=k/2,$ $l=k+j,$ $j\geq 0,$ and\newline
$\left(4\chi(\mathcal O_S)+t-K_S^2+8+2j-2n\right)k\leq 32\chi(\mathcal O_S)+4t-4K_S^2-8n,$\newline
with $n\leq j+7,$ where $n$ is the number of points of multiplicity $k/2.$
\item[e)] $r_m\leq k/2-2$ and\newline
$k\leq 5+\sqrt{1+8\chi(\mathcal O_S)},$ or\newline
$(4\chi(\mathcal O_S) + t - K_S^2)k\leq 32\chi(\mathcal O_S) + 4t - 4K_S^2.$
\end{description}
\item If $k\equiv 2\ ({\rm mod\ 4}),$ one of the following holds:
\begin{description}
\item[f)] $r_m=k/2+1$ and\newline
$(4\chi(\mathcal O_S)+t-K_S^2-2)k\leq 24\chi(\mathcal O_S) + 2t - 2K_S^2 - 20,$ with $t\geq 1,$\newline
or\newline
$(4\chi(\mathcal O_S) + t - K_S^2 - 8)k^2+(-32\chi(\mathcal O_S) - 4t + 4K_S^2 + 48)k+$\newline $80\chi(\mathcal O_S) + 4t - 4K_S^2 - 96\leq 0,$ with $t\geq 2;$
\item[g)] $r_m\leq k/2-1$ and\newline
$k\leq 5+\sqrt{1+8\chi(\mathcal O_S)},$ or\newline
$\left( 4\chi(\mathcal O_S)+t-K_S^2-6 \right)k\leq 24\chi(\mathcal O_S) + 2t - 2K_S^2 - 28.$
\end{description}
\end{description}
\end{proposition}
\noindent {\bf Proof:}
Let $H, G$ be as defined in Proposition \ref{rito} and let
$$P_1(l,r_m,G,H,k):=(2l-G)(k-r_m-2)-H,$$
$$P_2(l,r_m,G,H,k):=(2l-G)\big((r_m-4)(k-r_m)+(r_m-2)(k-r_m-2)\big)-H(2r_m-6).$$
From Lemma \ref{mainlemma}, $$P_1\leq 0\ \ \ {\rm and}\ \ \ P_2\leq 0.$$
\begin{description}
\item[a)] Let $n$ be the number of $(k/2+1,k/2+1)$ points. From Proposition \ref{semiXG}, b), d), $n=1$ or $2$. From Proposition \ref{rito}, we have
$$\sum (r_i-2)(k-r_i-2)=H'\ \ \ {\rm and}\ \ \ 2l=G'+\sum(r_i-2),$$
where
$$H'=H-n\left( k/2(k/2-4)+(k/2-2)^2 \right),\ \ \ G'=G+n(k-2)$$ and $r_i\leq k/2,$ $\forall i.$

The result follows from $$P_1(k+2,k/2,G',H',k)\leq 0.$$
\item[b)] From Proposition \ref{semiXG}, there are at most $(k/2+1,k/2+1)$ singularities. The inequality $$P_2(k+4,k/2+2,G,H,k)\leq 0$$ gives the result.
\item[c)] Let $n$ be the number of points of multiplicity $k/2$ and $m$ be the number of $(k/2-1,k/2-1)$ singularities.
From Proposition \ref{semiXG2}, c), $n=0$ or $1.$

If $n=0,$ then $r_m=k/2$ implies $m\geq 1$ (thus $t\geq 1$).
From $$P_2(k-2,k/2,G,H,k)\leq 0$$ one gets the first inequality.

Suppose $n=1.$
Notice that, as shown in the proof of Proposition \ref{semiXG2}, c), the point of multiplicity $k/2$ is obtained from the blow-up of
$\mathbb P^2$ at a point of type $(k/2-1,k/2-1).$ Hence $t\geq 1.$

Let $$H':=H-(k/2-2)^2,\ \ \ G'=G+k/2-2.$$

If $m=0,$ then $$P_1(k-2,k/2-2,G',H',k)\leq 0$$ implies the second inequality. 

If $m>0,$ then $$P_2(k-2,k/2,G',H',k)\leq 0$$ gives the third inequality.
In this case $t\geq 2.$

\item[d)] Let $j:=l-k$ and let $n$ be the number of points $x_i$ (possibly infinitely near) such that $r_i=k/2$.
From Proposition \ref{rito}, we have
$$\sum (r_i-2)(k-r_i-2)=H'\ \ \ {\rm and}\ \ \ 2l=G'+\sum(r_i-2),$$
where
$$H'=H-n(k/2-2)^2,\ \ \ G'=G+n(k/2-2)$$ and $r_i\leq k/2-2,$ $\forall i.$

The inequality $$P_1(k+j,k/2-2,G',H',k)\leq 0$$ gives $$\left(4\chi(\mathcal O_S)+t-K_S^2+8+2j-2n\right)k\leq 32\chi(\mathcal O_S)+4t-4K_S^2-8n.$$

It only remains to show that $n\leq j+7.$

One can verify, using the double cover formulas (see e.g. \cite{BHPV}), that $n\geq j+8$ implies $\chi(\mathcal O_S)<1,$ except for $n=8,$ $l=k$ and $n=10,$ $k=12,$ $l=14.$ 
We {\em claim} that in these cases $K_S^2\leq 0.$ This is impossible because $S$ is of general type.

{\em Proof of the claim}:\\
From the double cover formulas one gets that $\chi(\mc O_S)\leq 2$ and there is at least a $(-2)$-curve $A$ contained in $B$, otherwise $K_S^2\leq 0.$
One has $$B\equiv -\frac{k}{2}K_W+(l-k)\widetilde F+\sum\left(\frac{k}{2}-r_i\right)E_i,$$ where $\widetilde F$ is the total transform of $F$ and each $E_i$ is an exceptional divisor with self-intersection $-1.$
Since $AB=-2,$ $AK_W=0,$ $l\geq k$ and $r_i\leq k/2$ $\forall i,$ we have $AE_i<0$ for some $i$ such that $r_i<k/2.$ The only possibility is the existence of a $(3,3)$-point in $\b B$ and $\chi(\mc O_S)=1.$
But the imposition of such a singularity in the branch locus decreases the self-intersection of the canonical divisor by 1.
%
%
%
\item[e)] From Proposition \ref{semiXG2}, b), $l\geq k/2+r_m-2.$ Let $$f(r_m):=P_1(k/2+r_m-2,r_m,G,H,k).$$
We have $$f(r_m)=-2r_m^2+br_m+c\leq 0,$$ where $$b=4\chi(\mathcal O_S) + t - K_S^2 - k + 8$$ and $$c=k^2 - 10k -8\chi(\mathcal O_S) + 24.$$
Suppose that $c=f(0)>0$ (i.e. $k>5+\sqrt{1+8\chi(\mathcal O_S)}$). Then $f(r_m)$ has exactly one positive root $x.$ One has $$4x-b=\sqrt{b^2+8c}$$ and $k/2-2\geq r_m\geq x$ implies that
$$(4(k/2-2)-b)^2\geq b^2+8c.$$ This inequality gives the result.

\item[f)] Let $n$ be the number of points of type $(k/2,k/2).$

If $n=1,$ we proceed as in a).

If $n>1,$ the inequality is given by $$P_2(k,k/2+1,G,H,k)\leq 0.$$

\item[g)] It is analogous to the proof of e).\Qed
\end{description}

\section{Proof of main results}\label{PT}

\noindent{\bf Proof of Theorem \ref{thm1}:}\newline
Consider the parabola given by $f(x)=ax^2+bx+c,$ with $a>0.$ 
If $f(k)\leq 0,$ $f(z)\geq 0$ and $z\geq -b/2a$ (the first coordinate of the vertex), then $k\leq z.$

This fact and Proposition \ref{mainprop} imply that, if $K_S^2<4\chi(\mathcal O_S)-6,$ one of the following holds:

\begin{description}
\item[] \bf a) $k\leq\dfrac{16\chi(\mathcal O_S)-16}{4\chi(\mathcal O_S)-K_S^2-6}$\ \ \ \ \ \ \ \ \ \ \bf b) $k\leq\dfrac{16\chi(\mathcal O_S)}{4\chi(\mathcal O_S)+t-K_S^2-8},$ $t\geq 2$
\item[] \bf c) $k\leq 4+\dfrac{16\chi(\mathcal O_S)}{4\chi(\mathcal O_S)+t-K_S^2-4},$ $t\geq 1$
\item[] \bf c') $k\leq 4+\dfrac{16\chi(\mathcal O_S)-4}{4\chi(\mathcal O_S)+t-K_S^2-5},$ $t\geq 2$\ \ \ \ \bf d) $k\leq 4+\dfrac{16\chi(\mathcal O_S)-32}{4\chi(\mathcal O_S)-K_S^2-6}$
\item[] \bf e) $k\leq 5+\sqrt{1+8\chi(\mathcal O_S)}$\ \ \ \ \ \ \ \ \ \ \bf e') $k\leq 4+\dfrac{16\chi(\mathcal O_S)}{4\chi(\mathcal O_S)-K_S^2}$
\item[] \bf f) $k\leq 2+\dfrac{16\chi(\mathcal O_S)-16}{4\chi(\mathcal O_S)-K_S^2-1}$\ \ \ \ \ \bf f') $k\leq 2+\dfrac{16\chi(\mathcal O_S)-16}{4\chi(\mathcal O_S)+t-K_S^2-8},$ $t\geq 2$
\item[] \bf g) $k\leq 5+\sqrt{1+8\chi(\mathcal O_S)}$\ \ \ \ \ \ \ \ \ \bf g') $k\leq 2+\dfrac{16\chi(\mathcal O_S)-16}{4\chi(\mathcal O_S)-K_S^2}$
\end{description}

We want to show that $k$ is not greater than
\\\\
$\max \left\{ \frac{16\chi(\mathcal O_S)}{4\chi(\mathcal O_S)-K_S^2-6},4+\frac{16\chi(\mathcal O_S)-32}{4\chi(\mathcal O_S)-K_S^2-6} ,4+\frac{16\chi(\mathcal O_S)}{4\chi(\mathcal O_S)-K_S^2-3}, 5+\sqrt{1+8\chi(\mathcal O_S)} \right\}.$
\\\\
The result follows easily. Just notice that
$$4\chi(\mathcal O_S)-K_S^2-6\leq 8\Longrightarrow
2+\dfrac{16\chi(\mathcal O_S)-16}{4\chi(\mathcal O_S)-K_S^2-6} \leq 
\dfrac{16\chi(\mathcal O_S)}{4\chi(\mathcal O_S)-K_S^2-6}$$
and
$$4\chi(\mathcal O_S)-K_S^2-6\geq 8\Longrightarrow
2+\dfrac{16\chi(\mathcal O_S)-16}{4\chi(\mathcal O_S)-K_S^2-6} \leq 
4+\dfrac{16\chi(\mathcal O_S)-32}{4\chi(\mathcal O_S)-K_S^2-6}.$$\Qed
\noindent{\bf Proof of Proposition \ref{thm2}:}\newline
Let $(\alpha),$ $(\beta)$ be the equations of Proposition \ref{rito}, a), b), respectively. One has that $[(\alpha)+(k-10)(\beta)]/8$ is equivalent to
\begin{equation}\label{eqq1}
\frac{1}{8}\sum(r_i-2)(8-r_i)=15+K_S^2-t-3\chi(\mathcal O_S)-\frac{1}{4}(k-10)(l-10)
\end{equation}
and $(\beta)$+(\ref{eqq1}) is equivalent to
\begin{equation}\label{eqq2}
\chi(\mathcal O_S)=1+\frac{1}{4}(k-2)(l-2)-\frac{1}{8}\sum r_i(r_i-2).
\end{equation}
Now it suffices to show that $r_m\leq 8.$

Suppose that $K_S^2 < 3\chi(\mathcal O_S)-6.$

From \cite[Theorem 1]{Xi2} one gets that if $\chi(\mathcal O_S)\geq 54,$ then $S$ has a pencil of hyperelliptic curves of genus $\leq 6.$ In this case $k\leq 14,$ thus $r_m\leq k/2+2$ implies $r_m\leq 8.$

From the proof of Theorem \ref{thm1} we obtain that if $\chi(\mathcal O_S)\leq 31,$ then one of the possibilities below occur. In all cases $r_m\leq 8.$
\begin{description}
\item[] \bf a) b) $k<16,$ $r_m<8;$
\item[] \bf c) c') d) $k\leq 18,$ $r_m=k/2\leq 8;$
\item[] \bf e) $k\leq 20,$ $r_m\leq k/2-2\leq 8;$ \ \ \ \ \ \bf e') $k\leq 16,$ $r_m\leq k/2-2\leq 6;$
\item[] \bf f) $k\leq 14,$ $r_m=k/2+1\leq 8;$\ \ \ \ \ \ \bf f') $k\leq 16,$ $r_m=k/2+1\leq 8;$
\item[] \bf g) $k\leq 18,$ $r_m\leq k/2-1\leq 8;$\ \ \ \ \ \ \bf g') $k\leq 14,$ $r_m\leq k/2-1\leq 6.$
\end{description}

Suppose now that $32\leq\chi(\mathcal O_S)\leq 53.$
From Theorem \ref{thm1} we get that $k\leq 18$ or $k\leq 5+\sqrt{1+8\chi(\mathcal O_S)}.$ In this last case $r_m\leq k/2-1$ (see Proposition \ref{mainprop} e), g)).
Thus we have $r_m\leq 18/2+2$ or $r_m\leq 24/2-1.$ Since $r_m$ is even, $r_m\leq 10.$

Let $N_j$ be the number of points $x_i$ such that $r_i=j.$
We have $$\sum(r_i-2)\geq 8N_{10}+6N_8$$ and, from (\ref{eqq1}), $$8N_{10}\geq (k-10)(l-10)-32.$$
Using Proposition \ref{rito}, b) and the assumption $\chi(\mathcal O_S)\geq 32,$ we obtain $$2l+2k\geq 15+(k-10)(l-10)+6N_8.$$
This is equivalent to
\begin{equation}\label{eq2}
(k-12)(l-12)\leq 29-6N_8.
\end{equation}
Suppose $r_m=10.$ Then Propositions \ref{semiXG} and \ref{semiXG2} give two possibilities:
\begin{description}
\item[$\cdot$] $k=16,$ $l\geq k+2=18,$ there is a singularity of type $(9,9)$ ($N_8\geq 1$);
\item[$\cdot$] $k\geq 18,$ $l\geq k/2+r_m-2\geq 17.$
\end{description}
Both cases contradict (\ref{eq2}). We conclude that $r_m\leq 8.$\Qed
\noindent{\bf Proof of Theorem \ref{thm3}:}\newline
First we claim that if $A$ is a $(-2)$-curve contained in $B,$ the image $\b A$ of $A$ in $\m F_e$ does not intersect a negligible singularity of $\b B,$
unless $\b A$ is the negative section of $\m F_1$ and the only singularity of $\b B$ is a double point in $C_0$ (this corresponds to a smooth branch curve in $\m P^2$). 
In fact otherwise there is a $(-1)$-curve $E$ such that $AE=1$ or 2. If $AE=1,$ then $A+E$ can be contracted to a smooth point of the branch curve $\b B\subset\m F_e.$ This is impossible because the canonical resolution blows-up only singular points of $\b B.$ Suppose $AE=2.$ The inverse image of $A$ is a $(-1)$-curve which contracts to a smooth point of $S.$ The inverse image of $E$ is then contracted to a curve $\widehat E$ with arithmetic genus $1$ and $\widehat E^2=2.$ We obtain from the adjunction formula that $K_S\widehat E=-2,$ which is a contradiction because $S$ is of general type.

Recall that $t:=K_S^2-K_{S''}^2.$
The following holds:
\begin{enumerate}
\item[(1)] $l\geq k/2$\newline
	(Because $\b BC_0=l-ek/2\geq -e$ and $\b BC_0$ is even.);
\item[(2)] $l=k/2 \Longleftrightarrow (t=2 \wedge N_4=N_6=N_8=0)$\newline
	(In this case $e=1$ and $\b BC_0=0$.);
\item[(3)] $l=k/2+2 \Longrightarrow (N_6=N_8=0 \wedge t\geq N_4 \wedge (t=N_4 \vee N_4>1))$;\newline
	(If $N_4\ne 0,$ this corresponds to a branch curve in $\m P^2$ with $N_4$ points of type $(3,3)$ (see Proposition \ref{semiXG2}, c)).);
\item[(4)] $l=k-2 \wedge t=0 \Longrightarrow k/2$ even;\newline
	(As in $(1)$, $l\geq ek/2-e,$ hence $e\leq 2.$ We have $e=1$ because $t=0,$ thus $l$ even implies $k/2$ even.);
\item[(5)] $l<k-2 \Longrightarrow l-k/2$ even;\newline
	(As in $(1)$, $l\geq ek/2-e,$ thus $e=1$ and then $l-k/2=\b BC_0$ is even.)
\item[(6)] $t=1 \wedge N_4=N_6=N_8=0 \Longrightarrow l=k-2$.\newline
	(If there are only negligible singularities, $t=1$ is only possible if the negative section of $\m F_2$ is an isolated component of the branch locus.)
\end{enumerate}

For given values of $K_S^2-3\chi(\mc O_S)$ and $k,$ we want to choose the solution of the equation given in Proposition \ref{thm2}, b) which maximizes the value of $\chi(\mc O_S),$ given by the equation in Proposition \ref{thm2}, c). We can assume $N_6=N_8=0.$

It suffices to compute the numerical possibilities for Proposition \ref{thm2}, b), c) which satisfy conditions $(1),\ldots,(6)$.
We note the following:\\
since $k\geq 12,$ \cite[Thm. 1]{Xi2} implies $\chi(\mathcal O_S)\leq 69,$ then Theorem \ref{thm1} gives $k\leq 28$;\\
$l\geq k/2,$ $k\geq 12$ and (\ref{eqq1}) imply $-7\geq K_S^2-3\chi(\mathcal O_S)\geq -18+t+N_4,$ thus $K_S^2-3\chi(\mathcal O_S)\geq -18,$ $t\leq 11$ and $N_4\leq 11.$

A simple algorithm is available at\\ \verb+ http://home.utad.pt/~crito/magma_code.html +

The existence is easy to verify. All cases can be constructed as double covers of $\m P^2,$ $\m F_0,$ $\m F_1$ or $\m F_2.$
The table below contains information about $l$ or the degree of the branch curve in $\m P^2$ and about the singularities of the branch curve, if any.

\begin{table}[H]\scriptsize
\centering
\begin{tabular}{l l||l|l|l|l|l}
        ~               & g & ~                           & ~                          & ~                 & ~                           & ~ \\
        $K^2-3\chi$ & ~ &      -7                         & -8                         & -9                & -10                         & -11  \\\hline\hline
        ~               & 5 & $\m F_0,$ $l=26$            &  $\m F_0,$ $l=24$          & $\m F_0,$ $l=22$  &  $\m F_1,$ $l=20$           & $\m F_0,$ $l=18$ \\\hline
        ~               & 6 & $\m F_0,$ $l=18$            &  $\m F_1,$ $l=17$          & $\m F_0,$ $l=16$  &  $\m F_1,$ $l=15$           & $\m F_0,$ $l=14$ \\\hline
        ~               & 7 & $\m F_1,$ $l=14,$ $(3,3)$   &  $\m F_2,$ $l=14$          & $\m F_1,$ $l=14$  &  $\m F_1,$ $l=12,$ $(3,3)$  & $\m F_1,$ $l=12,$ $(4)$ \\\hline
        ~               & 8 & $\m F_1,$ $l=13,$ $(3,3)$   &  $\m F_1,$ $l=13,$ $(4)$   & $\m F_1,$ $l=13$  &  ~                          & $\m P^2,$ $20,$ $(3,3)$   \\\hline
        ~               & 9 & ~                           & $\m P^2,$ $22,$ $(3,3)$    & ~                 &  $\m F_1,$ $l=12$           & ~    \\\hline
        ~               & 10 & ~                          & ~                          & ~                 &  $\m P^2,$ $22$             & ~    \\     
\end{tabular}
\end{table}

\begin{table}[H]\scriptsize
\centering
\begin{tabular}{l l||l|l|l|l|l}
        ~                & g    & ~                  & ~                          & ~                  & ~                 & ~  \\
        $K^2-3\chi$      & ~    & -12                & -13                        & -14                & -15               & -16 \\\hline\hline
        ~                & 5    & $\m F_0,$ $l=16$   & $\m F_0,$ $l=14$           & $\m F_0,$ $l=12$   & $\m F_1,$ $l=10$  & $\m F_1,$ $l=8$  \\\hline
        ~                & 6    & $\m F_1,$ $l=13$   & $\m F_1,$ $l=11,$ $(4)$    & $\m F_1,$ $l=11$   & ~                 & $\m F_1,$ $l=9$  \\\hline
        ~                & 7    & $\m F_1,$ $l=12$   & $\m P^2,$ $18,$ $(3,3)$    & ~                  & $\m F_1,$ $l=10$  & $\m P^2,$ $16$  \\\hline
        ~                & 8    & ~                  & $\m F_1,$ $l=11$           & ~                  & $\m P^2,$ $18$    & ~   \\\hline
        ~                & 9    & ~                  & $\m P^2,$ $20$             & ~                  & ~                 & ~   \\\hline
        ~               & 10    & ~                  & ~                          & ~                  & ~                 & ~   \\     
\end{tabular}
\end{table}

\Qed

\bibliography{ReferencesRito}

\bigskip
\bigskip

\noindent Carlos Rito
\\ Departamento de Matem\' atica
\\ Universidade de Tr\' as-os-Montes e Alto Douro
\\ 5001-801 Vila Real
\\ Portugal
\\
{\it e-mail:} crito@utad.pt
\\\\

\noindent Mar\'ia Mart\'i S\'anchez

\noindent {\it e-mail:} mmartisanchez@educa.madrid.org

\end{document}